\documentclass{amsart}

 \theoremstyle{definition}
 
 \theoremstyle{remark}
 
\numberwithin{equation}{section}

 \newcommand{\D}{\mathfrak{D}}

 \newcommand{\DPRM}{\textsc{dprm }}
 
 \newcommand{\abs}[1]{\left\vert#1\right\vert}
 \newcommand{\set}[1]{\left\{#1\right\}}

\begin{document}
\title{Representations of $\Omega$ in number theory: finitude versus parity}
\author{Toby Ord}
\address{Department of Philosophy, The University of Melbourne, 3010, Australia}
\email{t.ord@pgrad.unimelb.edu.au}
\author{Tien D. Kieu}
\address{Centre for Atom Optics and Ultrafast Spectroscopy, Swinburne University of Technology,  Hawthorn 3122, Australia}
\email{kieu@swin.edu.au}


\begin{abstract}
We present a new method for expressing Chaitin's random real, $\Omega$, through Diophantine equations. Where Chaitin's method causes a particular quantity to express the bits of $\Omega$ by fluctuating between finite and infinite values, in our method this quantity is always finite and the bits of $\Omega$ are expressed in its fluctuations between odd and even values, allowing for some interesting developments. We then use exponential Diophantine equations to simplify this result and finally show how both methods can also be used to create polynomials which express the bits of $\Omega$ in the number of positive values they assume.
\end{abstract}

\maketitle
\begin{center}
\today
\end{center}



\section{Recursive Enumerability, Algorithmic Randomness and $\Omega$}
\label{section:Randomness}

One of the most startling recent developments in the theory of computation is the discovery of the number $\Omega$, through the subfield of algorithmic information theory. $\Omega$ is a real number between 0 and 1 which was introduced by G.~J.~Chaitin~\cite{Chaitin:1975} as an example of a number with two conflicting properties: it is both recursively enumerable and algorithmically random. Very roughly, this means that $\Omega$ has a simple definition and can be computed in the limit from below, yet we can determine only finitely many of its digits with certainty---for the rest we can do no better than random.

Understanding the full importance of these properties requires some familiarity with the recursive functions---commonly presented through models of computation such as Turing machines or the lambda calculus. For the purposes of algorithmic information theory, however, it is convenient to abstract some of the details from these models and consider a programming language in which the (partial) recursive functions are represented by finite binary strings.%
\footnote{For more details see Chaitin~\cite{Chaitin:1987}.}
These strings are just programs for a universal Turing machine (or universal lambda expression) and they take input in the form of a binary string then output another binary string or diverge (fail to halt). For convenience, we will often consider these inputs and outputs to encode tuples of positive integers.

On top of this simplified picture of computation, we impose one restriction which is necessary for the development of algorithmic information theory (and hence $\Omega$). The set of strings that encode the recursive functions must be prefix-free. This means that no program can be an extension of another, and thus each program is said to be self-delimiting. As algorithmic information theory is intricately linked with communication as well as computation, this is quite a natural constraint---if you wish to use a permanent binary communication channel, then you need to know when the end of a message has been reached and this cannot be done if some messages are extensions of others.

There are many prefix-free sets that one could choose and many recursive mappings between these and the recursive functions. These different choices of `programming language' lead to different values of $\Omega$, but this does not matter much as almost all of its significant properties will remain the same regardless. However, to allow talk of $\Omega$ as a specific real number we will use the same language as Chaitin~\cite{Chaitin:1987}.

Now that we have explained what we mean by a programming language, we can give a quick overview of computability in terms of programs. A program computes a set of $n$-tuples if, when provided with input $\langle x_1,\ldots,x_n \rangle$, it returns 1 if this is a member of the set and 0 otherwise. A program computes an infinite sequence if, when provided with input $n$, it returns the value of the $n$-th element in the sequence. A program computes a real, $r$, if it computes a sequence of rationals $\set{r_n}$ which converges to $r$ and $\abs{r - r_n} < \frac{1}{2^n}$. These sets, sequences and reals that are computed by programs are said to be recursive.

There are also many sets, sequences and reals that cannot be computed, but can be approximated in an important way. A program semi-computes a set of $n$-tuples if, when provided with input $\langle x_1,\ldots,x_n \rangle$, it returns 1 if this is a member of the set and diverges otherwise. A program semi-computes an infinite sequence of bits if, when provided with input $n$, it returns 1 if the $n$-th bit in the sequence is 1 and diverges otherwise. A program semi-computes a real, $r$, if, when provided with input $n$, it computes a rational number, $r_n$, where $\set{r_n}$ converges to $r$ from below. These sets, infinite bitstrings and reals that are semi-computed by programs are said to be recursively enumerable or r.e.

There is an important point that needs to be made concerning reals and their representations. Each real number between 0 and 1 has a binary expansion: a binary point followed by an infinite sequence of bits that represents the real.%
\footnote{For numbers that can be expressed with a representation ending in an infinite string of 0's, there is another representation ending in an infinite sequence of 1's, but we shall remove this ambiguity by only using representations with an infinite number of 0's. This will not affect the important reals in this paper, $\Omega$ and $\tau$, as they are irrational and thus have unique representations regardless.}
Throughout this paper, we shall be making considerable use of the binary expansions of real numbers so it is important to point out an oddity in the definitions above: a real is recursive if and only if its binary expansion is recursive, but a real may be r.e.~even if its binary expansion is not r.e. We shall thus take care to distinguish the weaker property of being an \emph{r.e.~real} from the stronger one of being a \emph{real whose binary expansion is r.e.}
\\

An example of a real that is r.e.~but not recursive is $\tau$: the real number between 0 and 1, whose $k$-th digit is 1 if the $k$-th program (in the usual lexical ordering of finite bitstrings) halts when given the empty string as input and 0 if the $k$-th program diverges. Equivalently:
\begin{equation}
\tau = \sum_{p_{n} \textrm{ halts}}2^{-n}
\end{equation}

$\tau$ is an r.e.~real because there is a computable sequence of rationals $\set{\tau_i}$, where
\begin{equation}
\tau_i = \sum_{\stackrel{\scriptstyle n\le i}{\scriptstyle p_n \textrm{ halts in } \le i \textrm{ steps}}} \!\!\!\!\!\!\!\!\!\!\!\!2^{-n}
\label{taui}
\end{equation} 
such that $\set{\tau_i}$ converges to $\tau$ from below.

Furthermore, it is clear that the binary representation of $\tau$ is also r.e.~because there is a program that simulates the $k$-th program, halting if and only if it does. This program is a slightly modified universal program that first determines the bits of the $k$-th program and then simulates it.

$\tau$ is not recursive, however, because if a program could compute it to arbitrary accuracy, it would determine whether each program halts or not when given the empty string as input. This is known as the \emph{blank tape problem} and is easily shown to be equivalent to the more general \emph{halting problem}---`does a given program halt on a given input?'. The halting problem is fundamental to the theory of computation and is the most famous problem that cannot be recursively solved. $\tau$ merely encodes the information necessary to solve the halting problem into the binary expansion of a real number and thus provides a very simple example of a non-computable real to which we can contrast the more exotic properties possessed by $\Omega$.

$\Omega$ encodes the halting problem in a more subtle way: it is the \emph{halting probability}. We could, theoretically, generate a random program one bit at a time, by flipping a fair coin and writing down a 1 when it comes up heads and a 0 for tails---stopping if we reach a valid program. The chance of generating any given $n$ bit program is therefore $\frac{1}{2^n}$. $\Omega$ is the chance that this method of random program construction generates a program that halts. Letting $\abs{p}$ represent the size of $p$ in bits, we can also express $\Omega$ as
\begin{equation} \label{eq:Omega}
\Omega = \sum_{p \textrm{ halts}}2^{-\abs{p}}
\end{equation}

As was the case for $\tau$, there is a computable sequence of rationals $\set{\Omega_i}$, where
\begin{equation}
\Omega_i = \sum_{\stackrel{\scriptstyle \abs{p}\le i}{\scriptstyle p \textrm{ halts in } \le i \textrm{ steps}}} \!\!\!\!\!\!\!\!\!\!\!\!2^{-\abs{p}}
\label{Omegai}
\end{equation}
which converges to $\Omega$ from below, showing it to be an r.e.~real. However, we shall see shortly that the binary representation of $\Omega$ is \emph{not} r.e.

A real is said to be \emph{algorithmically random}~\cite{Chaitin:1987} if and only if the `algorithmic complexity' of each $n$-bit initial segment of its binary expansion becomes and remains arbitrarily greater than $n$.%
\footnote{This is only one of four common definitions of algorithmic randomness, however, all have been shown to be equivalent.}
In other words a real, $r$, is algorithmically random if and only if any program that has access to outside advice in the form of binary messages requires more than $n$ bits of advice to compute the first $n$ bits of $r$'s binary expansion (for all values of $n$ above some threshold).%
\footnote{The reason that slightly more than $n$ bits of advice are needed is because in algorithmic information theory the advice comes in self-delimiting messages (which are actually programs that generate the advice---like self-extracting archives) and in order to be self-delimiting, these messages need slightly more bits than they would otherwise. In general, an $n$ bit string requires about $(n + \log{n})$ bits. Chaitin~\cite{Chaitin:1987} provides further details.}
Thus a random real is one for which only finitely many prefixes of its binary expansion can be compressed.

It is easy to see that a random real cannot have an r.e.~binary expansion. Let $x$ be an arbitrary real whose binary expansion is r.e. By definition, there must be a program, $p_x$, that takes a positive integer, $k$, and halts if and only if the $k$-th bit of $x$ is 1. To determine $n$ bits of $x$, we just need to know how many of these $n$ values of $k$ make $p_x$ halt. We could then simply run $p_x$ on all the values of $k$ and stop when this many have halted, knowing that no more will halt and thus determining the $n$ bits of $x$. Since all positive integers less than $n$ can be encoded in $\log{n}$ bits (rounding up), we only need to send a message of about $(\log{n} + \log{\log{n}})$ bits. In this manner, any prefix of $x$ can be significantly compressed, so $x$ cannot be random.

Because of this, we can see that $\tau$ too is not random. However, Chaitin~\cite{Chaitin:1987} has proven that $\Omega$ \emph{is} random and so cannot be compressed in this manner.%
\footnote{Indeed, it has since been shown through the work of R.~Solovay, C.~S.~Calude, P.~Hertling, B.~Khoussainov, Y.~Wang and T.~A.~Slaman that the only r.e.~random reals are $\Omega$'s for different programming languages. See Calude~\cite{Calude:2002} for more details.}
For sufficiently high values of $n$, $n$ bits of $\Omega$ provide $n$ bits of algorithmically incompressible information.
\\

In addition to recursive incompressibility, random reals are also characterised by \emph{recursive unpredictability}~\cite{Chaitin:1987}. Consider a `predictive' program that takes a finite initial segment of an infinite bitstring and returns a value indicating either `the next bit is 1', `the next bit is 0' or `no prediction'. If any such program is run on all finite prefixes of the binary expansion of a random real and makes an infinite amount of predictions, the limiting relative frequency of correct predictions approaches $\frac{1}{2}$. In other words when any program is used to predict infinitely many bits of a random real, such as $\Omega$, it does no better than random---even with information about all the prior bits.

The power of this unpredictability can be seen when compare the predictability of $\tau$. In this case, the predictive program can easily predict an infinite amount of bits with no errors. This is because infinitely many bits of $\tau$ are 'easy' to compute. For example, consider the halting behaviour of Turing machines: there are infinitely many Turing machines which have no loops in their transition graphs and thus cannot possibly diverge. When the predictive program is asked to predict the $n$-th bit of $\tau$, it can just check to see if the $n$-th program corresponds to such a machine, returning `the next bit is 1' if it does and `no prediction' otherwise.\footnote{From the definition of binary programs in algorithmic information theory, there must be a recursive mapping between programs and Turing machines (or any such model).}
\\

With its inherent incompressibility and unpredictability, $\Omega$ really does go beyond the type of uncomputability present in a more typical non-recursive real such as $\tau$. However, its contrasting property of being an r.e.~real makes $\Omega$ seem to be just beyond our reach. In the next section, we will introduce Diophantine equations and show how these can be used to bring uncomputability into the more classical field of number theory. Then, in Section~\ref{section:OmegaDiophantine}, we will show two ways of using Diophantine equations to bring $\Omega$ and randomness to number theory---Chaitin's original method and our new technique.


\section{Diophantine Equations and Hilbert's Tenth Problem}
\label{section:Diophantine}

A Diophantine equation is a polynomial equation in which all of the coefficients and variables take only positive integer values. Many natural phenomena with discrete quantities are modelled well by Diophantine equations and they occur frequently in number theory. It is often convenient to express a Diophantine equation with all terms on the left hand side:
\begin{equation}
D(x_{1},\ldots,x_{m}) = 0
\end{equation}
Here $D$ is a polynomial of $x_1,\ldots,x_m$ in which the coefficients can take both positive and negative integer values.

The number of solutions for a Diophantine equation varies widely. For example, $3x_1 + 6 = 0$ has one solution, while $x_1 x_2 - 2 = 0$ has two and $x_1 x_2 - x_2 = 0$ has infinitely many. Some however, such as $2 - 3x_1 = 0$, have no solutions at all. There are many different methods for deciding whether Diophantine equations of certain forms have solutions and determining what these solutions are, but there has been a great desire for a single method that takes an arbitrary Diophantine equation and determines whether or not it has solutions. In 1900, David Hilbert~\cite{Hilbert:1902} gave the problem of finding such a method as the tenth in his famous list of important problems to be addressed by mathematicians in the 20th Century. Since then, the task of finding this method has become known simply as Hilbert's Tenth Problem.

Another area of research concerns families of Diophantine equations. A family of Diophantine equations is a relation of the form:
\begin{equation}
D(a_{1},\ldots,a_{n},x_{1},\ldots,x_{m}) = 0
\end{equation}
in which we distinguish between two types of variable. The variables $x_{1},\ldots,x_{m}$ are called unknowns, while $a_{1},\ldots,a_{n}$ and called parameters. By assigning values to each of the parameters (and treating them as constants), we pick out an individual Diophantine equation from the family. For example, the family $a_1 - 3x_1 = 0$ consists of the equations: $1 - 3x_1 = 0$, $2 - 3x_1 = 0$, $3 - 3x_1 = 0$ and so on.

Each family of Diophantine equations is naturally associated with a certain set of $n$-tuples of positive integers, $\D$, in the following manner:
\begin{equation}
\langle a_{1},\ldots,a_{n} \rangle  \in \D
\quad \Longleftrightarrow \quad
\exists x_{1}\ldots{}x_{m}{D(a_{1},\ldots,a_{n},x_{1},\ldots,x_{m})=0}
\end{equation}
In other words, a tuple is in the set if the equation it corresponds to has a solution. Such sets are said to be Diophantine or to have a Diophantine representation. For example, the set of all multiples of 3 is Diophantine because it is represented by the family $a_1 - 3x_1 = 0$.

Over the 1950's and 1960's, M.~Davis, H.~Putnam and J.~Robinson established several important results regarding which sets are Diophantine. Their key result concerned a characterisation, not of Diophantine sets, but their close relation: \emph{exponential} Diophantine sets.

A family of exponential Diophantine equations is a relation of the form:
\begin{equation}
D(a_{1},\ldots,a_{n},x_{1},\ldots,x_{m},2^{x_{1}},\ldots,2^{x_{m}}) = 0
\end{equation}
where $D$ is once again a polynomial, but now some of its variables are exponential functions of others. Davis, Putnam and Robinson~\cite{DavisPutnamRobinson:1961} used this additional flexibility to show that all r.e.~sets are exponential Diophantine. It had long been known that all exponential (and standard) Diophantine sets are r.e.~because it is trivial to write a program that searches for a solution to a given equation and halts if and only if it finds one. Therefore, the new result meant that the exponential Diophantine sets were precisely the r.e.~sets.

In 1970, Yu.~Matiyasevich~\cite{Matiyasevich:1993} completed the final step, proving that all exponential Diophantine sets are also Diophantine and thus that the Diophantine sets are exactly the r.e.~sets---a result now known as the \DPRM Theorem.

The \DPRM Theorem provides an intimate link between Diophantine equations and computability, reducing the task of determining whether a set has a Diophantine representation to a matter of programming. For instance, there is a program that takes a single input $k$ and halts if and only if the $k$-th bit of $\tau$ is 1. Thus, the set of positive integers that includes $k$ if and only if the $k$-th program halts is an r.e.~set and via the \DPRM Theorem, there is a family of Diophantine equations with a parameter $k$, that has solutions if and only if the $k$-th program halts.

This family of equations provides an example of uncomputability in number theory and shows that Hilbert's Tenth Problem must be recursively undecidable because a program that finds whether arbitrary Diophantine equations have solutions could be used to determine the bits of $\tau$ and thus to solve the halting problem. Indeed, it was long known that the recursive undecidability of Hilbert's Tenth Problem would follow immediately from the \DPRM Theorem and this was the main motivation for its proof---the Diophantine representations for all other r.e.~sets being largely a bonus.


\section{Expressing Omega Through Diophantine Equations}
\label{section:OmegaDiophantine}

While the \DPRM Theorem demonstrates the existence of $\tau$ and uncomputability in number theory, it also denies the possibility of finding a similar family of Diophantine equations expressing $\Omega$ and randomness. This is due to the fact discussed in Section \ref{section:Randomness} that, while $\Omega$ is an r.e.~real, its sequence of bits is \emph{not} r.e. However, the \DPRM Theorem only prohibits a direct Diophantine representation of $\Omega$ and says nothing about the more subtle properties of Diophantine equations in which these bits could perhaps be encoded.

Chaitin~\cite{Chaitin:1987} takes such an approach. While there is no program of one variable, $k$, that halts if and only if the $k$-th bit of $\Omega$ is 1, Chaitin provides a program, $P$, that takes two variables, $k$ and $N$, and computes $\Omega$ somewhat less directly. For a given value of $k$, $P$ can be thought of as making an infinite series of `guesses' as to the value of the $k$-th bit of $\Omega$---when $P$ is run on $k$ and $N$, it gives the $N$-th guess as to the $k$-th bit of $\Omega$. What is impressive is that $P$ gets infinitely many of these guesses right and only finitely many wrong.

How does $P$ do this? It simply computes the sequence $\{\Omega_i\}$ discussed in Section~\ref{section:Randomness} until it gets to $\Omega_N$ and then returns the $k$-th bit of $\Omega_N$. Just as $\{\Omega_i\}$ forms a sequence of approximations to $\Omega$, so the $k$-th bit of each $\{\Omega_i\}$ forms a sequence of approximations to the $k$-th bit of $\Omega$.

Consider this $k$-th bit of each $\{\Omega_i\}$ as $i$ is increased. This bit could change between 0 and 1 many times, but since $\{\Omega_i\}$ approaches $\Omega$, it must eventually remain fixed, at which point it must have the same value as the $k$-th bit of $\Omega$. Therefore, if the $k$-th bit of $\Omega$ is 1, the $k$-th bit of $\{\Omega_i\}$ must be 0 for only finitely many values of $i$, and so $P$ must return 0 for finitely many values of $N$ and 1 for infinitely many. On the other hand, if the $k$-th bit of $\Omega$ is 0, then the $k$-th bit of $\{\Omega_i\}$ must be 1 for only a finite number of values of $i$ and $P$ must return 1 for finitely many values of $N$ and 0 for infinitely many. Either way, as $N$ increases, the output of $P$ applied to $k$ and $N$ limits to the $k$-th bit of $\Omega$.

It may seem as though this program is computing the bits of $\Omega$ but this is not quite the case. $P$ just computes the $N$-th `guess' of the $k$-th bit. From the infinite sequence of such guesses, the $k$-th bit could be determined but $P$ does not and cannot put the guesses together like that---it just returns one of them.

Since recursive functions are just a special type of r.e.~function, we can apply the \DPRM Theorem and see that there must be a family of Diophantine equations
\begin{equation}
\chi_1(k,N,x_1,\ldots,x_{m}) = 0
\end{equation}
that has solutions for given values of $k$ and $N$ if and only if $P$ returns 1 when provided with these as input. For a given value of $k$, there are solutions for infinitely many values of $N$ if and only if the $k$-th bit of $\Omega$ is 1.

Thus, by using a more subtle property of the family of Diophantine equations, Chaitin was able to show that algorithmic randomness occurs in number theory: as $k$ is varied, there is simply no recursive pattern to whether this family of equations has solutions for finitely or infinitely many values of $N$.
\\

By modifying Chaitin's method slightly, we can find a new way of expressing the bits of $\Omega$ through a family of Diophantine equations~\cite{OrdKieu:2003}. Consider a new program, $Q$, that also takes inputs $k$ and $N$, and begins to compute the sequence $\{\Omega_i\}$. For each value of $\Omega_i$, $Q$ checks to see if it is greater than $\frac{N}{2^k}$, halting if this is so, and continuing through the sequence otherwise. Since $\{\Omega_i\}$ approaches $\Omega$ from below, we can see that $\Omega_i > \frac{N}{2^k}$ implies that $\Omega > \frac{N}{2^k}$ and conversely, if $\Omega > \frac{N}{2^k}$ there must be some value of $i$ such that $\Omega_i > \frac{N}{2^k}$. Therefore, $Q$ will halt on $k$ and $N$ if and only if  $\Omega > \frac{N}{2^k}$. Alternatively, we could say that $Q$ recursively enumerates the pairs $\langle k, N \rangle$ such that $\Omega > \frac{N}{2^k}$.

Just as we could determine the $k$-th bit of $\Omega$ from the number of values of $N$ that make $P$ return 1, so we can determine it from the number of values of $N$ for which $Q$ halts. In what follows, we shall refer to these quantities as as $p_k$ and $q_k$ respectively. 

Unlike $p_k$, $q_k$ is always finite. Indeed, an upper bound is easily found. Since $\Omega < 1$, only values of $k$ and $N$ such that $\frac{N}{2^k} < 1$ can possibly be less than $\Omega$ and thus make $Q$ halt. Since both $k$ and $N$ take only values from the positive integers we also know that $\frac{N}{2^k} > 0$ and thus for a given $k$, there are less than $2^{k}$ values of $N$ for which $Q$ halts and $q_k \in \set{0,1,\ldots,2^k-1}$.

From the value of $q_k$, it is quite easy to derive the first $k$ bits of $\Omega$. Firstly, note that $q_k$ is equal to the largest value of $N$ such that $\frac{N}{2^k} < \Omega$---unless there is no such $N$, in which case it equals 0. Either way, its value can be used to provide a very tight bound on the value of $\Omega$: $\frac{q_k}{2^k} < \Omega \le \frac{q_k + 1}{2^k}$. Since $\Omega$ is irrational, we can strengthen this to $\frac{q_k}{2^k} < \Omega < \frac{q_k + 1}{2^k}$, which means that the first $k$ bits of $\frac{q_k}{2^k}$ are exactly the first $k$ bits of $\Omega$.

This gives some nice results connecting $q_k$ and $\Omega$. The first $k$ bits of $\frac{q_k}{2^k}$ are just the bits of $q_k$ when written with enough leading zeros to make $k$ digits in total. Thus $q_k$, when written in this manner, provides the first $k$ bits of $\Omega$. Additionally, we can see that $q_k$ is odd if and only if the $k$-th bit of $\Omega$ is 1.

Now that we know the power and flexibility of $q_k$, it is a simple matter to follow Chaitin in bringing these results to number theory. The function computed by $Q$ is r.e.~so, by the \DPRM Theorem, there must be a family of Diophantine equations
\begin{equation}
\chi_2(k,N,x_1,\ldots,x_{m}) = 0
\end{equation}
that has a solution for specified values of $k$ and $N$ if and only if $Q$ halts when given these values as inputs. Therefore, for a particular value of $k$, this equation only has solutions for values of $N$ between 0 and $2^k - 1$ with the number of solutions, $q_k$, being odd if and only if the $k$-th bit of $\Omega$ is 1.

This new family of Diophantine equations improves upon the original one in a couple of ways. Whereas the first method expressed the bits of $\Omega$ in the fluctuations between a finite and infinite amount of values of $N$ that give solutions, the second keeps this value finite and bounded, with the bits of $\Omega$ expressed through the more mundane property of parity. It is the fact that this quantity is always finite that leads to many of the new features of this family of Diophantine equations. $p_k$ is infinite when the $k$-th bit of $\Omega$ is 1 and, since there is only one way in which it can be infinite, it can provide no more than this one bit of information. On the other hand, $q_k$ can be odd (or even) in $2^{k-1}$ ways, which is enough to give $k-1$ additional bits of information, allowing the first $k$ bits of $\Omega$ to be determined.

The fact that $q_k$ is always finite also provides a direct reduction of the problem of determining the bits of $\Omega$ to Hilbert's Tenth Problem. To find the first $k$ bits of $\Omega$, one need only determine for how many values of $N$ the new family of Diophantine equations has solutions. Since we know that there can be no solutions for values of $N$ greater than or equal to $2^k$, we could determine the first $k$ bits of $\Omega$ from the solutions to $2^k$ instances of Hilbert's Tenth Problem. In fact, we can lower this number by taking advantage of the fact that if there is a solution for a given value of $N$ then there are solutions for all lower values. All we need is to find the highest value of $N$ for which there is a solution and we can do this with a bisection search, requiring the solution of only $k$ instances of Hilbert's Tenth Problem.\footnote{For details see \cite{OrdKieu:2003}.}

Finally, the fact that $q_k$ is always finite allows the generalisation of these results from binary to any other base, $b$. If we replace all above references to $2^k$ with $b^k$ we get a new program, $Q_b$, with its associated family of Diophantine equations. For this family, the value of $q_k$ now gives us the first $k$ digits of the base $b$ expansion of $\Omega$: it is simply the base $b$ representation of $q_k$ with enough leading zeroes to give $k$ digits. The value of the $k$-th digit of $\Omega$ is simply $q_k$ mod $b$.
\\

Chaitin~\cite{Chaitin:1987} did not stop with his Diophantine representation of $\Omega$, but instead moved to exponential Diophantine equations where his result could be presented more clearly. He made this move to take advantage of the theorem that all r.e.~sets have \emph{singlefold} exponential Diophantine representations, where a representation is singlefold if each equation in the family has at most one solution.

We can denote the singlefold family of exponential Diophantine equations for the program $P$ by
\begin{equation}
\chi^e_1(k,N,x_1,\ldots,x_{m'}) = 0
\end{equation}
For a given $k$, this equation will have exactly one solution for each of infinitely many values of $N$ if the $k$-th bit of $\Omega$ is 1 and exactly one solution for each of finitely many values of $N$ if the $k$-th bit of $\Omega$ is 0. We can make use of this to express the bits of $\Omega$ through a more intuitive property.

If we treat $N$ in this equation as an unknown instead of a parameter, we get a new (very similar) family of exponential Diophantine equations with only one parameter
\begin{equation}
\label{equation:FinitudeExpDiophantine}
\chi^e_1(k,x_0,x_1,\ldots,x_{m'}) = 0
\end{equation}
Since the previous family was singlefold and $N$ has become another unknown, there will be exactly one solution to this single parameter family for each value of $N$ that gave a solution to the double parameter family. Thus, (\ref{equation:FinitudeExpDiophantine}) has infinitely many solutions if and only if the $k$-th bit of $\Omega$ is 1.

This same approach can be used with our method~\cite{OrdKieu:2003}. There is a two-parameter singlefold family of exponential Diophantine equations for $Q$ and this can be converted to a single parameter family of exponential Diophantine equations
\begin{equation}
\chi^e_2(k,x_0,x_1,\ldots,x_{m'}) = 0
\end{equation}
with between 0 and $2^{k} - 1$ solutions, the quantity being odd if and only if the $k$-th bit of $\Omega$ is 1.
\\

Finally, we have also shown~\cite{OrdKieu:2003} that both Chaitin's finitude-based method and our parity-based method can be used to generate polynomials for $\Omega$. For a given family of Diophantine equations with two parameters, 
\begin{equation}
\label{equation:GenericDiophantine}
D(k,N,x_1,\ldots,x_m) = 0
\end{equation}
we can construct a polynomial, $W$, where
\begin{equation}
W(k,x_{0},x_{1},\ldots,x_{m})
\equiv
x_{0}\left(1-(D(k,x_{0},x_{1},\ldots,x_{m}))^{2}\right).
\end{equation}
Note that the parameter, $N$, is again treated as an unknown and thus denoted $x_0$.

If we restrict the values of the variables to positive integers then, for a given $k$, this polynomial takes on exactly the set of all values of $N$ for which (\ref{equation:GenericDiophantine}) has solutions. We can thus use this method on $\chi_1 = 0$ and $\chi_2 = 0$, generating polynomials that express $p_k$ and $q_k$ in the number of distinct positive integer values they take on for different values of $k$. We therefore have a polynomial whose number of distinct positive integer values fluctuates from odd to even and back in an algorithmically random manner as a parameter $k$ is increased.


\section*{Acknowledgements}
We would like to thank Gregory Chaitin for helpful discussion and Cristian Calude for kindly offering to present this paper for us at DMTCS'03. TDK wishes to acknowledge the continuing support of Peter Hannaford.

\bibliography{diophantine,ait,hypercomputation}
\bibliographystyle{plain}

\end{document}